\documentclass[12pt]{article}
\usepackage{amssymb}

\newtheorem{definition}{Definition}[section]
\newtheorem{theorem}[definition]{Theorem}
\newtheorem{lemma}[definition]{Lemma}

\newtheorem{note}[definition]{Note}

\def\C{\mathbb C}
\def\K{\mathbb K}

\def\Z{\mathbb Z}

\newcommand{\beast}{\begin{eqnarray*}}
\newcommand{\eeast}{\end{eqnarray*}}

\begin{document}
\newenvironment{proof}{\noindent{\it Proof\/}:}{\par\noindent $\Box$\par}

\title{ \bf Two linear transformations each tridiagonal \\ 
 \bf with respect to an eigenbasis of the other;\\
 an overview 
}
\author{Paul Terwilliger}
\date{}
\maketitle
\begin{abstract} 
Let $\K$ denote a field 
 and let $V$ denote a  
vector space over $\K$ with finite positive dimension.
We consider an ordered pair of linear transformations 
$A:V\rightarrow V$ and $A^*:V\rightarrow V$ 
that 
 satisfy conditions (i), (ii) below. 
\begin{enumerate}
\item There exists a basis for $V$ with respect to which
the matrix representing $A$ is irreducible tridiagonal and
the matrix representing $A^*$ is diagonal.
\item There exists a basis for $V$ with respect to which
the matrix 
representing $A$ is diagonal and 
the matrix representing $A^*$ is irreducible tridiagonal.
\end{enumerate}
We call such a pair a {\it Leonard pair} on $V$.
We give an overview of the theory of Leonard pairs.
\end{abstract}

\section{Leonard pairs}
\medskip

\medskip
\noindent We begin by recalling the notion of a Leonard pair.
We will use the following terms. Let $X$ denote a square matrix.
Then $X$ is called {\it tridiagonal} whenever each nonzero
entry lies on either the diagonal, the subdiagonal, or the
superdiagonal. Assume $X$ is tridiagonal. Then $X$ is
called {\it irreducible} whenever each entry on
the subdiagonal is nonzero and each entry on the superdiagonal
is nonzero.

\medskip
\noindent We now define a Leonard pair.  For the rest of this paper
$\K$ will denote a field.

\begin{definition} 
\label{def:lprecall}
\cite{LS99}
\rm
Let  
 $V$ denote a  
vector space over $\K$ with finite positive dimension.
By a {\it Leonard pair} on $V$
we mean an ordered pair of linear transformations 
$A:V\rightarrow V$ and $A^*:V\rightarrow V$ 
that 
 satisfy conditions (i), (ii) below. 
\begin{enumerate}
\item There exists a basis for $V$ with respect to which
the matrix representing $A$ is irreducible tridiagonal and
the matrix representing $A^*$ is diagonal.
\item There exists a basis for $V$ with respect to which
the matrix 
representing $A$ is diagonal and 
the matrix representing $A^*$ is irreducible tridiagonal.
\end{enumerate}
\end{definition}

\begin{note}
\rm
According to a common notational convention,
$A^*$ denotes the conjugate transpose of $A$.
We are not using this convention.
In a Leonard pair $A,A^*$ the linear transformations $A,A^*$
are arbitrary subject to (i), (ii) above.
\end{note}

\noindent 
\begin{note}
\rm
Our use of the name ``Leonard pair'' is motivated by a connection
to a theorem of D. Leonard 
\cite[p.  260]{BanIto}, \cite{Leon}
involving the $q$-Racah and related
polynomials of the Askey Scheme.
\end{note}


\section{An example of a Leonard pair} 
\medskip
\noindent
Here is an example of a Leonard pair.
Set 
$V={\K}^4$ (column vectors), set 
\beast
A = 
\left(
\begin{array}{ c c c c }
0 & 3  &  0    & 0  \\
1 & 0  &  2   &  0    \\
0  & 2  & 0   & 1 \\
0  & 0  & 3  & 0 \\
\end{array}
\right), \qquad  
A^* = 
\left(
\begin{array}{ c c c c }
3 & 0  &  0    & 0  \\
0 & 1  &  0   &  0    \\
0  & 0  & -1   & 0 \\
0  & 0  & 0  & -3 \\
\end{array}
\right),
\eeast
and view $A$ and $A^*$  as linear transformations on $V$.
We assume 
the characteristic of $\K$ is not 2 or 3 to ensure
$A$ is irreducible.
Then the pair $A, A^*$ is a Leonard
pair on $V$. 
Indeed 
condition (i) in Definition
\ref{def:lprecall}
is satisfied by the basis for $V$
consisting of the columns of the 4 by 4 identity matrix.
To verify condition (ii), we display an invertible  matrix  
$P$ such that 
$P^{-1}AP$ is 
diagonal and 
$P^{-1}A^*P$ is
irreducible tridiagonal.
Set 
\beast
P = 
\left(
\begin{array}{ c c c c}
1 & 3  &  3    &  1 \\
1 & 1  &  -1    &  -1\\
1  & -1  & -1  & 1  \\
1  & -3  & 3  & -1 \\
\end{array}
\right).
\eeast
 By matrix multiplication $P^2=8I$, where $I$ denotes the identity,   
so $P^{-1}$ exists. Also by matrix multiplication,    
\begin{equation}
AP = PA^*.
\label{eq:apeq}
\end{equation}
Apparently
$P^{-1}AP$ is equal to $A^*$ and is therefore diagonal.
By (\ref{eq:apeq}) and since $P^{-1}$ is
a scalar multiple of $P$, we find
$P^{-1}A^*P$ is equal to $A$ and is therefore irreducible tridiagonal.  Now 
condition (ii) of  Definition 
\ref{def:lprecall}
is satisfied
by the basis for $V$ consisting of the columns of $P$. 

\section{Leonard systems}

When working with a Leonard pair, it is often convenient to
consider a closely related and somewhat more abstract concept
called a {\it Leonard system}.
In order to define this we recall a few terms.
Let $V$ denote a vector space over $\K$ with finite positive
dimension. Let $\mbox{End}(V)$ denote the $\K$-algebra
consisting of the linear transformations from $V$ to $V$.
For $A \in 
 \mbox{End}(V)$, by the {\it eigenvalues} of $A$ we mean 
 the roots of the characteristic polynomial of $A$. These eigenvalues
 are contained in the algebraic closure of $\K$.
 We say $A$ is {\it multiplicity-free}
 whenever the eigenvalues of $A$ are mutually distinct and contained in
 $\K$. Assume for the moment that $A$ is multiplicity-free.
Let $\theta_0,\theta_1, \ldots, \theta_d$ denote an ordering
of the eigenvalues of $A$. For $0 \leq i \leq d$ let 
$v_i$ denote a nonzero vector in  $V$ that is an eigenvector
for $A$ with eigenvalue $\theta_i$.
Observe the sequence $v_0, v_1, \ldots, v_d$ is a basis for $V$.
For $0 \leq i \leq d$ define 
$E_i \in \mbox{End}(V)$
so that $E_iv_j=\delta_{ij}v_j$ for $0 \leq j\leq d$.
We call $E_i$ the {\it primitive idempotent} of $A$ associated with
$\theta_i$.
\begin{definition}
\label{def:ls}
\cite{LS99}
\rm
Let $d$ denote a nonnegative integer and 
let $V$ denote a vector space over $\K$ with dimension $d+1$.
By a {\it Leonard system} on $V$ we mean a sequence
$\Phi=(A;A^*; $ $\lbrace E_i \rbrace_{i=0}^d; 
\lbrace E^*_i \rbrace_{i=0}^d)$ that satisfies 
conditions (i)--(v) below.
\begin{enumerate}
\item Each of $A,A^*$ is a multiplicity-free element of
 $\mbox{End}(V)$.
\item $E_0, E_1, \ldots, E_d$ is an ordering of the primitive idempotents
of $A$.
\item $E^*_0, E^*_1, \ldots, E^*_d$ is an ordering of the primitive idempotents
of $A^*$.
\item 
${\displaystyle{
E^*_iAE^*_j = \cases{0, &if $\;|i-j|> 1$;\cr
\not=0, &if $\;|i-j|=1$\cr}
\qquad \qquad 
(0 \leq i,j\leq d).
}}$
\item
${\displaystyle{
 E_iA^*E_j = \cases{0, &if $\;|i-j|> 1$;\cr
\not=0, &if $\;|i-j|=1$\cr}
\qquad \qquad 
(0 \leq i,j\leq d).
}}$
\end{enumerate}
We call $d$ the {\it diameter} of $\Phi$. 
\end{definition}
Leonard pairs and Leonard systems are related as follows.

\begin{theorem}
\cite{TLT:split}
For $A,A^*$ in 
 $\mbox{End}(V)$,
 the pair $A,A^*$
is a Leonard pair on $V$ if and only if
the following (i), (ii) hold.
\begin{enumerate}
\item Each of $A,A^*$ is multiplicity-free.
\item There exists an ordering $E_0, E_1, \ldots, E_d$
of the primitive idempotents of $A$
and there exists an ordering
 $E^*_0, E^*_1, \ldots, E^*_d$
of the primitive idempotents of $A^*$ such that
$(A;A^*;$ $\lbrace E_i \rbrace_{i=0}^d; 
\lbrace E^*_i \rbrace_{i=0}^d)$ is a Leonard system on $V$.
\end{enumerate}
\end{theorem}

\section{Leonard pairs from 24 points of view}

\medskip
Let $A, A^*$ denote a Leonard 
pair on $V$.
We describe  
24 bases for $V$ on which $A,A^*$ act in an attractive fashion.

\medskip
\noindent 
To describe the bases
we will use the following terms.
Let $X$ denote a square matrix. Then $X$ is called
{\it lower bidiagonal} 
 whenever each nonzero
entry lies on either the diagonal or the subdiagonal.
We say $X$ is 
{\it upper bidiagonal} 
 whenever the transpose of $X$ is lower bidiagonal.
Let $v_0, v_1, \ldots, v_d$ denote
a  basis for $V$.
By 
the {\it inversion} of 
this basis
we mean the basis $v_d, v_{d-1}, \ldots, v_0$.

\medskip
\noindent The 24 bases are described by the diagram below.
In that diagram each  vertex represents one  of the
24 bases.
The shading on the vertex  indicates the nature
of the
$A,A^*$ action.
\begin{enumerate}
\item
Black: 
 $A$ is 
 diagonal
and  $A^*$
is irreducible tridiagonal.
\item Green: 
$A$ is lower bidiagonal and $A^*$ is upper bidiagonal.
\item
Red: 
$A$ is upper bidiagonal and $A^*$ is lower bidiagonal.
\item Yellow:
 $A$ is irreducible tridiagonal 
and $A^*$ is diagonal.
\end{enumerate}

\noindent
For each pair  
 of bases in the 
diagram that are connected by an arc,
consider the transition matrix from
one of these bases to the other. 
The shading on the arc indicates the nature of
this transition matrix. 
\begin{enumerate}
\item Solid arc: Transition matrix
is diagonal. 
\item
Dashed arc:
Transition matrix
is lower triangular.
\item
Dotted arc: The bases are the inversion
of one another.
\end{enumerate}

\bigskip
\input psfig.sty
\centerline{\psfig{figure=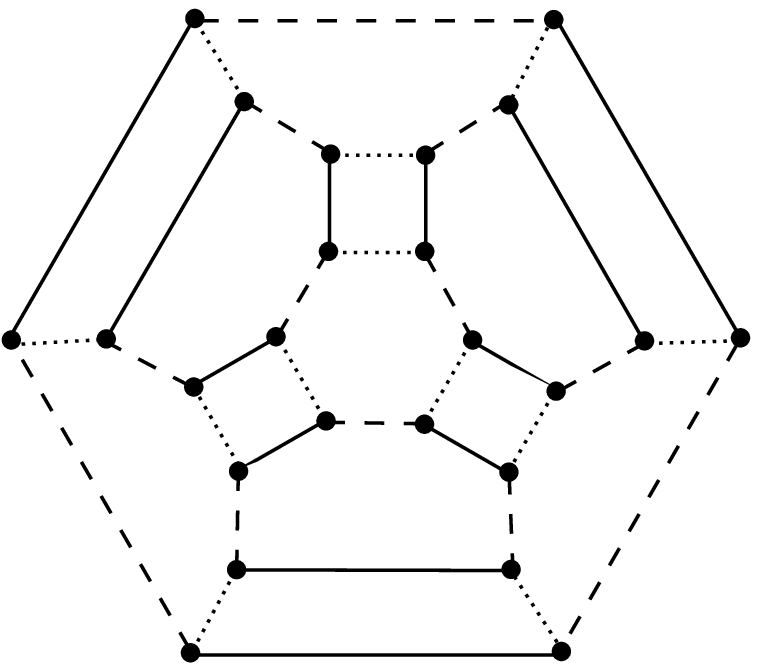,height=4.5cm}} 
\vskip .5cm

\bigskip
\noindent See
   \cite{LS24} for more information concerning the 24 bases.

\section{The classifying space}

We will shortly give a classification of the Leonard systems.
In order to describe the result we recall the notion
of a {\it parameter array}.

\begin{definition} 
\label{def:pa}
\rm
Let $d$ denote a nonnegative integer.
By a {\it parameter array over $\K$ of diameter $d$}
we mean a sequence of scalars 
$
(\theta_i, \theta^*_i,i=0..d; \varphi_j, \phi_j,j=1..d)
$
taken from $\K$ which satisfy the following conditions
(PA1)--(PA5).
\begin{description}
\item[(PA1)]  
$
\theta_i \not=\theta_j, \qquad 
\theta^*_i \not=\theta^*_j \qquad \mbox{if}\quad i\not=j, \qquad \qquad
(0 \leq i,j\leq d)
$.
\\
\item[(PA2)]
$
\varphi_i \not=0, \qquad \phi_i \not=0 \qquad \qquad (1 \leq i \leq d).
$
\\
\item[(PA3)]
$\varphi_i = \phi_1 \sum_{h=0}^{i-1} 
\frac{\theta_h-\theta_{d-h}}{\theta_0-\theta_d}
+ (\theta^*_i-\theta^*_0)(\theta_{i-1}-\theta_d)
\qquad \qquad (1 \leq i \leq d)$.
\\
\item[(PA4)] 
$\phi_i = \varphi_1 \sum_{h=0}^{i-1} 
\frac{\theta_h-\theta_{d-h}}{\theta_0-\theta_d}
+ (\theta^*_i-\theta^*_0)(\theta_{d-i+1}-\theta_0)
\qquad \qquad (1 \leq i \leq d)$.
\\
\item[(PA5)]
The expressions
\begin{eqnarray*} 
\frac{\theta_{i-2}-\theta_{i+1}}{\theta_{i-1}-\theta_i},
\qquad 
\frac{\theta^*_{i-2}-\theta^*_{i+1}}{\theta^*_{i-1}-\theta^*_i}
\label{eq:betaplusone}
\end{eqnarray*}
are equal and independent of $i$ for $2 \leq i \leq d-1$.
\end{description}
\end{definition}

\section{The classification of Leonard systems}

In this section we give a bijection from
the set of parameter arrays to the set of isomorphism
classes of Leonard systems.
Let
$(\theta_i, \theta^*_i,i=0..d;\varphi_j,\phi_j,j=1..d)$
denote a parameter array over $\K$.
Define
\begin{eqnarray*}
A=\left(
\begin{array}{c c c c c c}
\theta_0 & & & & & {\mathbf 0} \\
1 & \theta_{1} &  & & & \\
& 1 & \theta_{2} &  & & \\
& & \cdot & \cdot &  &  \\
& & & \cdot & \cdot &  \\
{\mathbf 0}& & & & 1 & \theta_d
\end{array}
\right),
\quad A^*=\left(
\begin{array}{c c c c c c}
\theta^*_0 &\varphi_1 & & & & {\mathbf 0} \\
 & \theta^*_1 & \varphi_2 & & & \\
&  & \theta^*_2 & \cdot & & \\
& &  & \cdot & \cdot &  \\
& & &  & \cdot & \varphi_d \\
{\mathbf 0}& & & &  & \theta^*_d
\end{array}
\right).
\end{eqnarray*}
Observe 
$A$ (resp. $A^*$) is multiplicity free
with eigenvalues
$\theta_0, \theta_1, \ldots, \theta_d$
(resp.
$\theta^*_0, \theta^*_1, \ldots, \theta^*_d$).
For $0 \leq i \leq d$ let $E_i$ (resp. $E^*_i$)
denote the primitive idempotent of $A$ (resp. $A^*$)
associated with $\theta_i$ (resp. $\theta^*_i$).
Then the sequence
$(A;A^*; $ $\lbrace E_i \rbrace_{i=0}^d; 
\lbrace E^*_i \rbrace_{i=0}^d)$ 
is a Leonard system on $\K^{d+1}$.
This construction induces a bijection from
the set of parameter arrays to the set of isomorphism
classes of Leonard systems \cite{LS99}.

\section{Leonard pairs $A,A^*$ with $A$ lower bidiagonal and
$A^*$ upper bidiagonal}

\noindent 
Let $A, A^*$ denote
matrices in $\hbox{Mat}_{d+1}(\K)$. Let us assume 
$A$ is lower bidiagonal 
and $A^*$ is upper  bidiagonal. 
We give a necessary and sufficient condition
for $A,A^*$ to be  a Leonard pair.

\begin{theorem}
   \cite{conform}
\label{thm:lug}
Let $d$ denote a nonnegative integer and let $A, A^*$ denote
matrices in $\hbox{Mat}_{d+1}(\K)$. Assume $A$  
lower bidiagonal and $A^*$ is upper bidiagonal.
Then the following (i), (ii) are equivalent.
\begin{enumerate}
\item
The pair $A,A^*$ is a Leonard pair.
\item
There exists a parameter array
$(\theta_i, \theta^*_i, i=0..d;  \varphi_j, \phi_j, j=1..d)$
over $\K$ 
such that
\begin{eqnarray*}
&&A_{ii} =\theta_i,
\qquad \qquad 
A^*_{ii} =\theta^*_i \qquad \qquad (0 \leq i \leq d),
\label{eq:comb1}
\\
&&\qquad A_{i,i-1}A^*_{i-1,i} = \varphi_i \qquad \qquad (1 \leq i \leq d).
\label{eq:comb2}
\end{eqnarray*}
\end{enumerate}
\end{theorem}

\section{Leonard pairs $A,A^*$ with $A$ tridiagonal and $A^*$ diagonal}

\noindent 
Let $A, A^*$ denote
matrices in $\hbox{Mat}_{d+1}(\K)$. Let us assume 
$A$ is tridiagonal 
and $A^*$ is diagonal. 
We give  a necessary and sufficient condition for
 $A,A^*$ to be a Leonard pair. 

\begin{theorem}
   \cite{conform}
\label{thm:tdcrit}
Let $d$ denote a nonnegative integer and let $A, A^*$ denote
matrices in $\hbox{Mat}_{d+1}(\K)$. Assume $A$ is
tridiagonal and $A^*$ is diagonal.
Then the following (i), (ii) are equivalent.
\begin{enumerate}
\item
The pair $A,A^*$ is a Leonard pair.
\item
There exists
a parameter array 
$(\theta_i, \theta^*_i, i=0..d;  \varphi_j, \phi_j, j=1..d)$
over $\K$
such that
\begin{eqnarray*}
A_{ii} &=&
\theta_i + \frac{\varphi_i}{\theta^*_i-\theta^*_{i-1}}
+
 \frac{\varphi_{i+1}}{\theta^*_i-\theta^*_{i+1}}
\qquad \qquad  
(0 \leq i \leq d), 
\label{eq:aiiform}
\\
A_{i,i-1}A_{i-1,i}&=&
\varphi_i\phi_i \frac{\prod_{h=0}^{i-2}(\theta^*_{i-1}-\theta^*_h)
}
{\prod_{h=0}^{i-1}(\theta^*_{i}-\theta^*_h)
}
\, \frac{\prod_{h=i+1}^d (\theta^*_i-\theta^*_h)
}
{\prod_{h=i}^d(\theta^*_{i-1}-\theta^*_h)
}
\qquad  (1 \leq i \leq d),
\label{eq:crossprod}
\\
A^*_{ii} &=&\theta^*_i \qquad \qquad (0 \leq i \leq d).
\label{eq:ths}
\end{eqnarray*}
\end{enumerate}
\end{theorem}

\section{A characterization of the parameter arrays I}

In this section we characterize the parameter
arrays in terms of bidiagonal matrices.
We will refer to the following set-up.

\begin{definition}
\label{def:setup}
Let $d$ denote a nonnegative integer and let
$(\theta_i, \theta^*_i,i=0..d; \varphi_j, \phi_j,j=1..d)
$
denote a sequence of  scalars taken from $\K$. We assume
this sequence satisfies 
PA1 and PA2.
\end{definition}

\begin{theorem}
 \cite{TLT:array}
\label{thm:gmat}
With reference to Definition \ref{def:setup},
 the following (i), (ii) are equivalent.
\begin{enumerate}
\item
The sequence
$(\theta_i, \theta^*_i,i=0..d;\varphi_j,\phi_j,j=1..d)$
satisfies PA3--PA5.
\item
There exists an invertible matrix 
$G \in \mbox{Mat}_{d+1}(\K)$ such that both
\begin{eqnarray*}
G^{-1}
\left(
\begin{array}{c c c c c c}
\theta_0 & & & & & {\mathbf 0} \\
1 & \theta_{1} &  & & & \\
& 1 & \theta_{2} &  & & \\
& & \cdot & \cdot &  &  \\
& & & \cdot & \cdot &  \\
{\mathbf 0}& & & & 1 & \theta_d
\end{array}
\right) G
&=&\left(
\begin{array}{c c c c c c}
\theta_d & & & & & {\mathbf 0} \\
1 & \theta_{d-1} &  & & & \\
& 1 & \theta_{d-2} &  & & \\
& & \cdot & \cdot &  &  \\
& & & \cdot & \cdot &  \\
{\mathbf 0}& & & & 1 & \theta_0
\end{array}
\right),
\label{eq:g1}
\\
G^{-1}
\left(
\begin{array}{c c c c c c}
\theta^*_0 &\varphi_1 & & & & {\mathbf 0} \\
 & \theta^*_1 & \varphi_2 & & & \\
&  & \theta^*_2 & \cdot & & \\
& &  & \cdot & \cdot &  \\
& & &  & \cdot & \varphi_d \\
{\mathbf 0}& & & &  & \theta^*_d
\end{array}
\right) G
&=&
\left(
\begin{array}{c c c c c c}
\theta^*_0 &\phi_1 & & & & {\mathbf 0} \\
 & \theta^*_1 & \phi_2 & & & \\
&  & \theta^*_2 & \cdot & & \\
& &  & \cdot & \cdot &  \\
& & &  & \cdot & \phi_d \\
{\mathbf 0}& & & &  & \theta^*_d
\end{array}
\right).
\label{eq:g2}
\end{eqnarray*}
\end{enumerate}
\end{theorem}

\section{A characterization of the parameter arrays II}

\noindent In this section we characterize the parameter
arrays in terms of polynomials.
We will use the following notation.
Let $\lambda $ denote an indeterminate, and let
$\K \lbrack \lambda \rbrack $ denote the $\K$-algebra
consisting of all polynomials in $\lambda $  which have
coefficients in $\K$. From now on
all polynomials which we discuss
are assumed to lie in $\K\lbrack \lambda \rbrack $.

\begin{theorem}
 \cite{TLT:array}
\label{eq:mth}
With reference to Definition \ref{def:setup},
the following (i), (ii) are equivalent.
\begin{enumerate}
\item
The sequence
 $(\theta_i, \theta^*_i,i=0..d; \varphi_j, \phi_j,j=1..d)$
satisfies PA3--PA5.
\item For $0 \leq i \leq d$ the 
polynomial 
\begin{equation}
\sum_{n=0}^i \frac{
(\lambda-\theta_0)
(\lambda-\theta_1) \cdots
(\lambda-\theta_{n-1})
(\theta^*_i-\theta^*_0)
(\theta^*_i-\theta^*_1) \cdots
(\theta^*_i-\theta^*_{n-1})
}
{\varphi_1\varphi_2\cdots \varphi_n}
\label{eq:poly1}
\end{equation}
is a scalar multiple of the polynomial 
\begin{eqnarray*}
\sum_{n=0}^i \frac{
(\lambda-\theta_d)
(\lambda-\theta_{d-1}) \cdots
(\lambda-\theta_{d-n+1})
(\theta^*_i-\theta^*_0)
(\theta^*_i-\theta^*_1) \cdots
(\theta^*_i-\theta^*_{n-1})
}
{\phi_1\phi_2\cdots \phi_n}.
\label{eq:poly2}
\end{eqnarray*}
\end{enumerate}
\end{theorem}

\section{Some orthogonal polynomials of the Askey scheme}

There is a natural correspondence between 
Leonard systems and a class of orthogonal polynomials 
consisting of the $q$-Racah polynomials and some related
polynomials of the Askey scheme. This correspondence
is described as follows
 \cite{TLT:array}.

\medskip
\noindent
Let 
$\Phi$
denote a Leonard system.  
By the {\it polynomials which correspond to $\Phi$} we mean 
the polynomials 
from  
(\ref{eq:poly1}), where the parameter array from
that line is the one associated with 
$\Phi$.
The polynomials which correspond to a Leonard system
are listed below:

\medskip
\noindent
$q$-Racah, $q$-Hahn, dual $q$-Hahn,
$q$-Krawtchouk,
dual $q$-Krawtchouk,
quantum $q$-Krawtchouk, 
affine $q$-Krawtchouk, 
Racah, Hahn, dual-Hahn, Krawtchouk, Bannai/Ito, and  orphan polynomials.

\medskip
\noindent
The Bannai/Ito polynomials can be obtained from the $q$-Racah polynomials by letting $q$ 
tend to $-1$. The orphan polynomials have maximal degree 3 and 
exist for $\mbox{char}(\K)=2$  only. 

\medskip
\noindent See 
 \cite{KoeSwa} for information on the Askey scheme.

\section{The Askey-Wilson relations}

We turn our attention to the representation theoretic aspects
of Leonard pairs.

\begin{theorem}
\cite{aw}
\label{lptheorem}
Let $V$ denote a vector space over $\K$ with finite positive
dimension. Let $A,A^*$ denote
a Leonard pair on $V$. Then there
exists a sequence of scalars
$\beta,\gamma,\gamma^*,\varrho,\varrho^*$, $\omega,\eta,\eta^*$
taken from $\K$ such that both
\begin{eqnarray*}  \label{askwil1}
A^2 A^*-\beta A A^*\!A+A^*\!A^2-\gamma\left( A A^*\!+\!A^*\!A
\right)-\varrho\,A^* &=& \gamma^*\!A^2+\omega A+\eta\,I,\\
\label{askwil2} A^*{}^2\!A-\beta A^*\!AA^*\!+AA^*{}^2-
\gamma^*\!\left(A^*\!A\!+\!A A^*\right)-\varrho^*\!A &=&
\gamma A^*{}^2+\omega A^*\!+\eta^*I.
\end{eqnarray*}
The sequence is uniquely determined by the pair $A,A^*$ provided
the dimension of $V$ is at least $4$.
\end{theorem}

\noindent 
The following theorem is a kind of converse to Theorem
 \ref{lptheorem}.

\begin{theorem}
\cite{aw}
Let $V$ denote a vector space over $\K$ with finite positive
dimension. Let $A:V\to V$ and $A^*:V\to V$ denote linear
transformations. Suppose that:
\begin{itemize}
\item There exists a sequence of scalars
$\beta,\gamma,\gamma^*,\varrho,\varrho^*,\omega,\eta,\eta^*$ taken
from $\K$ which satisfies {\rm (\ref{askwil1})}, {\rm
(\ref{askwil2})}.
\item $q$ is not a root of unity, where $q+q^{-1}=\beta$.
\item Each of $A$ and $A^*$ is multiplicity-free.
\item There does not exist a subspace $W\subseteq V$ such that
$W\not=0$,
$W\not=V$,
$AW\subseteq W$,
$A^*W\subseteq W$.
\end{itemize}
Then $A,A^*$ is a Leonard pair on $V$.
\end{theorem}

\section{Leonard pairs and the Lie algebra $sl_2$}

\medskip
\noindent
In this section we assume
the field $\K$ is algebraically closed 
with characteristic zero.

\medskip
\noindent We recall the Lie algebra $sl_2=sl_2(\K)$. 
This algebra has a basis  $e,f,h$ satisfying
\beast
&&\lbrack h,e \rbrack  = 2e, \qquad  
\lbrack h,f \rbrack  = -2f, 
\qquad
\lbrack e,f \rbrack  = h,
\eeast
where $\lbrack\,,\, \rbrack $ denotes the Lie bracket.

\medskip
\noindent 
We recall the irreducible finite dimensional modules for $sl_2$.
\begin{lemma} \cite[p. 102]{Kassel} 
\label{lemvd}
There exists a family 
\begin{equation}
\label{sl2modlist}
V_d  \qquad \qquad d = 0,1,2\ldots
\end{equation}
of irreducible finite dimensional $sl_2$-modules with the following
properties. The module $V_d$ has a basis $v_0, v_1, \ldots, v_d$
satisfying $h v_i =(d-2i)v_i$ for $0 \leq i \leq d$,
$fv_i = (i+1)v_{i+1}$ for $0 \leq i \leq d-1$, $fv_d=0$, 
$ev_i = (d-i+1)v_{i-1}$ for $1 \leq i \leq d$,  $ev_0=0$.
Every irreducible finite dimensional $sl_2$-module
is isomorphic to exactly one of the modules 
in line (\ref{sl2modlist}).
\end{lemma}

\begin{theorem} 
\label{ex:slgen}
\cite[Ex. 1.5]{TD00}
Let $A$ and $A^*$ denote semi-simple 
elements in $sl_2$ and assume $sl_2$ is generated by
these elements. 
Let $V$ denote an irreducible finite dimensional
 module for $sl_2$. Then
the pair $A, A^*$ acts on $V$ as a Leonard  pair.
\end{theorem}

\noindent We remark the Leonard pairs in Theorem 
\ref{ex:slgen} correspond to the  Krawtchouk
polynomials
 \cite{KoeSwa}.

\section{Leonard pairs and $U_q(sl_2)$}

\medskip
\noindent In this section we assume $\K$ is algebraically closed.
We fix a nonzero scalar $q\in \K$ which is not a root of unity. 
We recall the quantum algebra
 $U_q(sl_2)$.

\begin{definition} 
\cite[p.122]{Kassel}
\label{def:uqsl2}
Let $U_{q}(sl_2)$ denote the associative $\K$-algebra 
with 1 generated by symbols $ e,  f, k, k^{-1}$ subject
to the relations
\beast
kk^{-1} = k^{-1}k= 1,
\eeast
\beast
 ke = q^2 ek,\qquad\qquad  kf = q^{-2}fk,
\eeast
\beast
ef - fe  = {{k-k^{-1}}\over {q - q^{-1}}}.
\eeast

\end{definition}
\noindent We recall the irreducible finite dimensional modules
for 
 $U_{q}(sl_2)$. We use the following notation.
\beast
\lbrack n \rbrack_q = 
{{q^n - q^{-n}}\over {q-q^{-1}}} \qquad 
\qquad n \in \Z .
\eeast

\begin{lemma} \cite[p. 128]{Kassel}
\label{lem:uqmods}
With reference to Definition \ref{def:uqsl2},
there exists a family 
\begin{eqnarray}
V_{\varepsilon,d} \qquad \quad 
\varepsilon \in \lbrace 1,-1\rbrace, \qquad \quad  d = 0,1,2\ldots
\label{eq:uqmods}
\end{eqnarray}
of 
 irreducible finite dimensional 
 $U_{q}(sl_2)$-modules
 with the following properties.
The module 
$V_{\varepsilon,d}$ has a basis $u_0, u_1, \ldots, u_d$
satisfying $ku_i=\varepsilon q^{d-2i}u_i$ for
$0 \leq i \leq d$, $fu_i=\lbrack i+1\rbrack_q u_{i+1}$
for $0 \leq i \leq d-1$, $fu_d=0$, 
 $eu_i=\varepsilon\lbrack d-i+1\rbrack_q u_{i-1}$
for $1 \leq i \leq d$, $eu_0=0$. 
Every irreducible finite dimensional
 $U_{q}(sl_2)$-module is isomorphic to exactly one of  
the modules $V_{\varepsilon,d}$.
(Referring to line
(\ref{eq:uqmods}), if 
$\K$ has characteristic 2
we interpret the set $\lbrace 1,-1 \rbrace $ as having 
a single element.) 
\end{lemma}


\begin{theorem}
\cite{cite38}, \cite{cite39},
\cite{Terint}
\label{ex:uqex}
Referring to Definition
\ref{def:uqsl2} and Lemma 
\ref{lem:uqmods},
let  $\alpha, \beta $ denote nonzero scalars in  $\K $ and 
 define $A$, $A^*$ as follows.
\begin{eqnarray*}
A= \alpha f + {{ k}\over {q-q^{-1}}},
\qquad \qquad 
A^*=\beta e + {{k^{-1}}\over {q-q^{-1}}}.
\label{eq:abdef}
\end{eqnarray*}
Let $d$ denote a nonnegative integer and choose 
$\varepsilon \in \lbrace 1, -1 \rbrace $.
 Then the pair $A, A^*$ acts on $V_{\varepsilon,d}$ as a Leonard pair
provided $\varepsilon \alpha \beta $ is not among
$q^{d-1}, q^{d-3},\ldots, q^{1-d}$. 
\end{theorem}

\noindent 
We remark  
the Leonard pairs in 
Theorem
\ref{ex:uqex}
correspond to the quantum $q$-Krawtchouk
polynomials \cite{KoeSwa}, 
\cite{cite37}. 

\section{Leonard pairs in  combinatorics}

Leonard pairs arise in many branches of combinatorics.
For instance they arise in
the theory of partially ordered sets (posets).
We illustrate this with a poset called the  
subspace lattice $L_n(q)$.

\medskip
\noindent 
In this section we assume our field $\K$ is the field $\C$ of
complex numbers.

\medskip
\noindent To define the subspace lattice
 we introduce a second field. 
Let $GF(q)$ denote a finite field of order $q$.
Let $n$ denote a positive integer
and let $W$ denote an $n$-dimensional vector space over $GF(q)$.
Let $P$ denote the set consisting of all subspaces of $W$.
The set $P$, together with the containment relation, is a poset
called $L_n(q)$. 

\medskip
\noindent
Using $L_n(q)$ we obtain a family of Leonard pairs as follows.
Let  $\C P$  denote
the vector space over $\C$ consisting of all formal $\C$-linear combinations
of elements of $P$. 
We observe $P$ is a basis for  
 $\C P$ so
the dimension of $\C P$ is equal to the cardinality of P.

\medskip
\noindent
We define  three linear transformations on $\C P$. We call these
$K$, $R$ (for ``raising''), $L$ (for ``lowering'').

\medskip
\noindent
We begin with $K$. For all $x \in P$,
\beast  
	  K x =   q^{{n/2} - dim\, x} x.
\eeast
Apparently each element of $P$ is an eigenvector for $K$.

\medskip
\noindent
To define $R$ and $L$ we use the following notation. 
For $x, y \in P$ we say $y$ {\it covers}  $x$
whenever (i) $x \subseteq y$ and  (ii) $\hbox{dim}\, y= 1+\hbox{dim}\,x$.

\medskip
\noindent The maps 
$R$ and $L$ are defined as follows. For all $x \in  P$, 
\beast
           R x = \sum_{y \;covers \; x} y.
\eeast  
Similarly
 \beast 
           L x =q^{(1-n)/2} \sum_{x \; covers \; y} y.
\eeast
(The scalar $q^{(1-n)/2}$ is included for aesthetic reasons.)

\medskip
\noindent We consider the properties of $K, R, L$.
From the construction we find $K^{-1}$ exists.
By combinatorial counting we  verify 
\beast
     &&     K L = q L K,   \qquad \qquad \qquad   K R = q^{-1} R K,
     \\
       && \qquad          L R - R L = \frac{K - K^{-1}}{q^{1/2} - q^{-1/2}}.
\eeast
We recognize these equations. They are the defining relations for 
$U_{q^{1/2}}(sl_2)$.  Apparently $K$, $R$, $L$ turn  $\C P$ into
a module for $U_{q^{1/2}}(sl_2)$. 

\medskip
\noindent
We now see how to get Leonard pairs from $L_n(q)$.
Let  $\alpha, \beta$ denote nonzero complex scalars and 
define $A$, $A^*$ as follows.
\beast
                  A =  \alpha R    + \frac{K}{q^{1/2}-q^{-1/2}},
        \qquad \qquad  A^* =  \beta  L    + \frac{K^{-1}}{q^{1/2}-q^{-1/2}}.
\eeast
To avoid 
degenerate situations we assume  
$\alpha \beta $  is not among $q^{(n-1)/2}, q^{(n-3)/2},\ldots, q^{(1-n)/2}$.

\medskip
\noindent The  $U_{q^{1/2}}(sl_2)$-module $\C P$ is 
completely reducible \cite[p. 144]{Kassel}.
In other words $\C P$ 
is a direct sum of irreducible
$U_{q^{1/2}}(sl_2)$-modules. On
each irreducible module in this sum
the pair $A, A^*$ acts as a Leonard pair. This follows from
Theorem 
\ref{ex:uqex}.

\medskip
\noindent
We just saw how the subspace lattice gives Leonard pairs. 
It is implicit in \cite{uniform} that
the following posets give 
Leonard pairs in a similar fashion:
the subset lattice,
the Hamming semi-lattice,
the attenuated spaces,
and the classical polar spaces. 
Definitions of these posets can be found in 
 \cite{uniform}.

\section{Further reading}

We mention some additional topics which are related 
to Leonard pairs.

\medskip
\noindent 
Earlier in this paper we obtained 
Leonard pairs from the irreducible finite dimensional modules
for the Lie algebra $sl_2$ and the quantum algebra $U_q(sl_2)$.
We cite some other 
algebras whose modules are related to Leonard pairs.
These are the Askey-Wilson algebra
\cite{GYZnature},
 \cite{GYLZmut},
\cite{GYZTwisted},
\cite{GYZlinear},
\cite{GYZspherical},
\cite{Zhidd}, 
\cite{ZheCart},
\cite{Zhidden}, the Onsager algebra
\cite{CKOnsn},
\cite{DateRoan2},
\cite{Dav},
\cite{Dolgra},
and the Tridiagonal algebra
\cite{TD00},
\cite{LS99},
\cite{qSerre}.

\medskip
\noindent 
We discussed how
certain classical posets give Leonard pairs.
Another combinatorial object which gives Leonard pairs is
 a  $P$- and $Q$-polynomial
association scheme \cite{BanIto}, \cite{bcn}, 
 \cite{TersubI}.
Leonard pairs have been used to describe
certain irreducible modules for the subconstituent algebra of
 these schemes
\cite{Cau}, \cite{CurNom}, \cite{Curspin}, 
\cite{TD00},
 \cite{TersubI}.

\medskip
\noindent 
The topic of Leonard pairs is closely related to the  work of 
Gr\"unbaum and Haine on the 
 ``bispectral problem''
\cite{GH7},
\cite{GH6}.
See
\cite{GH4},
\cite{GH5},
\cite{GH1}, 
\cite{GH3},
\cite{GH2} 
for related work.

\noindent Paul Terwilliger, Department of Mathematics, University of
Wisconsin, 480 Lincoln Drive, Madison, Wisconsin, 53706, USA \hfil\break
Email: terwilli@math.wisc.edu \hfil\break

\end{document}